\documentclass[11pt,a4paper]{article}
\usepackage{pifont}
\usepackage{bbding}

\usepackage{amssymb}
\usepackage{mathrsfs}
\usepackage{amsfonts}
\usepackage{amsthm,amscd,amsmath}
\usepackage{eufrak}
\usepackage{youngtab}
\usepackage{graphicx}
\usepackage[all]{xy}

\makeatletter
\def\@biblabel#1{#1}
\makeatother

\newtheorem{theorem}{Theorem}[section]
\newtheorem{lemma}[theorem]{Lemma}

\theoremstyle{definition}
\newtheorem{definition}[theorem]{Definition}

\theoremstyle{proposition}
\newtheorem{proposition}[theorem]{Proposition}

\theoremstyle{remark}
\newtheorem{remark}[theorem]{Remark}

\theoremstyle{corollary}
\newtheorem{corollary}[theorem]{Corollary}
\date{}

\begin{document}

\title{Canonical bases and quantum coordinate ring}

\author{Bin Li and Hechun Zhang\\
\small{Department of Mathematical Science}, \\
\small{Tsinghua University, Beijing, P.R. China, 100084}\\
\scriptsize{Email: hzhang@math.tsinghua.edu.cn,
 libin07@mails.tsinghua.edu.cn}}

\maketitle

\begin{abstract}
Some filtrations  of  the tensor product of a highest weight module
and a lowest weight module over quantum group $U_q(\mathfrak g)$ are
constructed in \cite{LZ:2009} and one can use them to define some ideals
of the modified quantized enveloping algebra. It is shown that the quotient algebras inherit canonical
bases from the modified quantized enveloping algebra and are dual to
the quantum coordinate ring defined by Kashiwara for symmetrizable
Kac-Moody algebra $\mathfrak g$.
\end{abstract}


\section{Introduction}

Quantum coordinate ring, or quantum function algebra, is the
q-deformed version of the coordinate ring associated to a Lie group
$G$. It can be viewed, in some sense, as an algebra dual to the
quantized enveloping algebra $U=U_q(\mathfrak g)$ and thus it is
natural to study  its structure and representations as well as its
$\Bbb Z$-form.\\
\indent There are various ways to define quantum coordinate ring
$C$. For any Kac-Moody algebra $\mathfrak g$ with a symmetrizable
Generalized Cartan matrix, M. Kashiwara defined in \cite{Kashiwara2:1993}
$C$ as the algebra generated by
 all coordinate functions of the $U$-modules in the category
$\mathcal O_{int}$ and moreover, there is an analogue of Peter-Weyl
theorem
\[C\cong \bigoplus_{\lambda\in P^+}V(\lambda)\otimes V(\lambda)^{\circ}\]
where $V(\lambda)$ is the irreducible integrable highest
weight $U$-module with highest weight $\lambda$ and
$V(\lambda)^{\circ}$ is its graded dual. In particular for
 $\mathfrak g$ of finite type, Lusztig  gave another equivalent definition \cite {Lusztig:1992}
 of the quantum coordinate ring
 and it is known that the quantum coordinate ring
  is exactly the Hopf
dual of $U$ in this case. In the present paper, we will follow
Lusztig's approach to
 define a quantum coordinate ring through
modified quantized enveloping algebra $\widetilde{U}$. \\
\indent Recall that G. Lusztig constructed in \cite{Lusztig:1992} a
canonical basis for the tensor product $V(\lambda)\otimes V(-\mu)$
as well as for $\widetilde{U}$. When $\mathfrak g$ is of finite
 type, he considered the subspace $\widetilde{U}^{\circ}$ of
the dual space of $\widetilde{U}$ spanned by the dual basis of the
canonical basis of $\widetilde{U}$ in \cite{Lusztig:2009}.
 The multiplication in $\widetilde{U}^{\circ}$ is defined through
  the coproduct on $\widetilde{U}$ to make it become an associative
   algebra which is proved later to be isomorphic to the quantum
    coordinate ring. In this way, the
    $\Bbb Z$-form of this ring is naturally defined \cite{Lusztig:2009}.
    But unfortunately this method is not valid for $\mathfrak g$ of other
     types and hence we need to do something more.\\
\indent In \cite{Lusztig:1992} Lusztig conjectured that
  for $\mathfrak g$ of finite type, there is a composition series
   of $V(\lambda)\otimes V(-\mu)$ compatible with the canonical
   basis. In \cite{Lusztig:1993}, Lusztig gave an inductive method to construct the composition series of any integrable module in category
   $\mathcal O_{int}$. A different approach
 to construct the composition series is
given in \cite{LZ:2009} based on the theory of crystal basis. With this
 method we can also construct a nice
filtration of $V(\lambda)\otimes V(-\mu)$ for $\mathfrak g$ of any
type such that the quotient of any two neighbors is either zero or an irreducible
integrable highest weight module. Using these filtrations, we define a subspace
$\widetilde{U}^{\prime}$ of $\widetilde{U}$ spanned by all canonical
 base elements $G(b)$ such that $b\in\widetilde{B}$ is contained in
a connected component not isomorphic to a highest weight crystal. It
is proved that $\widetilde{U}^{\prime}$ is a two-sided ideal of
$\widetilde{U}$. The quotient $\mathfrak{U}\triangleq
\widetilde{U}/\widetilde{U}^{\prime}$ is an associative algebra
which inherits from $\widetilde{U}$ a canonical basis. Let $\mathfrak{U}$
take the place of $\widetilde{U}$ and then we define, similar to
what  Lusztig did in \cite{Lusztig:2009}, an algebra which is
proved to be isomorphic to the quantum coordinate ring.\\
\indent The quantum coordinate ring considered in this paper
involves only integrable representations in category $\mathcal
O_{int}$. Thus it is exactly the algebra of strongly regular
functions on symmetrizable Kac-Moody group \cite{Kac:1983} when $q=1$. It
would be interesting if $\mathcal O_{int}$ is replaced by  some
larger categories. Namely, if there are more generators besides
coordinate functions of highest weight modules, say, those of lowest
weight modules, the generated subalgebra of $U^{\ast}$ will be much
more interesting. For $\mathfrak g$ of affine type, the structure of
level zero part of $\widetilde{U}$  was studied by Beck and Nakajima
\cite{BN:2004}. The authors of the present paper believe that their work is
helpful  to understand this coordinate ring of
affine type though this is not included in the
present paper.\\
\indent The paper is organized as follows. In section 2, we recall
 some definitions and facts about the crystal and canonical
  bases of $\widetilde{U}$. In particular, the action of Cartan involution
  on canonical basis of $\widetilde{U}$
  is studied through the bilinear form on $\widetilde{U}$. In section 3 and 4, some nice
  filtrations
   of the tensor product $V(\lambda)\otimes V(-\mu)$ are constructed. We then define $\mathfrak{U}$ to be the quotient of
     $\widetilde{U}$ and investigate its cell modules as in \cite{Lusztig:1993}. In the last section, an algebra
      dual to $\mathfrak{U}$ is defined and proved to be isomorphic to quantum coordinate
      ring.

\section{Preliminaries}

\subsection{Modified quantized enveloping algebra $\widetilde{U}$}

We denote by $\mathfrak g=\mathfrak g(A)$ any symmetrizable
Kac-Moody algebra of rank $n$. The set of simple roots is indexed by
 $I=\{1,\cdots,n\}$. Let $Q$ be its root lattice, i.e.
$Q=\bigoplus_{i\in I}\Bbb Z\alpha_i\subset \mathfrak h$
where $\mathfrak h$ is the Cartan subalgebra and $\alpha_i$ are the
simple roots. Let $\Pi^{\vee}=\{h_i\in\mathfrak h |\ i\in I\}$ be
the set of simple coroots. We choose $d_j\in\mathfrak h$,
$1\leqslant j\leqslant n-rank(A)$ such that $\Pi^{\vee}\bigcup
\{d_j\in\mathfrak h\ |\ 1\leqslant j\leqslant n-rank(A)\}$ forms a
basis of $\mathfrak h$. Set \[P^{\vee}=\bigoplus_{i\in I}\Bbb
Zh_i\bigoplus\bigoplus_{1\leqslant j\leqslant n-rank(A)}\Bbb Zd_j\subset
\mathfrak h.\] The weight lattice $P$ is defined as
$P=\{\lambda\in\mathfrak h^{\ast}\ |\ \lambda(h)\in\Bbb Z\ \forall h\in P^{\vee}\}.$
Let $Q^+$ and $P^+$ be the positive root lattice and the set of
dominant weights respectively. We define $P_0$ to be the subset of
$P^+$ consisting of weights $\mu$ such that
$\mu(h_i)=0\ \textrm{for\ all}\ i\in I.$ Let $W$ be the Weyl group
associated to $\mathfrak g$. There is a $W$-invariant symmetric
bilinear form $(\ ,\ )$ on $P\times P$ such that
\[2\frac{(\lambda,\alpha_i)}{(\alpha_i,\alpha_i)}=\lambda(h_i).\]
 Let $U_q(\mathfrak g)$ be the quantized enveloping algebra generated
 over $k=\Bbb Q(q)$ by $E_i$, $F_i$ and $q^{h}$
 for $i\in I$, $h\in P^{\vee}$ \cite{Kashiwara:1994}, which is denoted also by $U$ for simplicity. The subalgebras $U^+$, $U^0$ and
 $U^-$ are defined in the same way as in \cite{Kashiwara:1994}.
 For $\xi=\sum_{i\in I}n_i\alpha_i\in Q$, define the height of $\xi$ to be
$\sum_{i\in I}|n_i|$, denoted by $ht(\xi)$. Set $U_{\xi}=\{u\in U\ |\ q^{h}uq^{-h}=q^{\xi(h)}u\}.$
The filtration $F=(F_n)_{n\in\Bbb Z_+}$ of $U^{\pm}$ is
defined by
\[F_n(U^{\pm})=\bigoplus_{ht(\xi)\leqslant n}U^{\pm}_{\xi}.\]
Denote by  $\widetilde{U}_q(\mathfrak g)$ or
 simply
 $\widetilde{U}$ the modified quantized enveloping algebra
 generated by $U_q(\mathfrak g)a_{\lambda}$ for $\lambda\in P$ subject to the relations:
\[q^{h}a_{\lambda}=q^{\lambda(h)}a_{\lambda},\ a_{\lambda}a_{\mu}=\delta_{\lambda,\mu}a_{\lambda},\
 ua_{\lambda}=a_{\lambda+\xi}u\ \textrm{for}\ u\in U_{\xi}.\]
 Note that $\widetilde{U}=\bigoplus_{\lambda\in P}Ua_{\lambda}.$

There is an anti-automorphism (resp. automorphism)  of $U$, denoted by $\ast$
(resp. $\omega$), such that
\[E_i^{\ast}=E_i,\ F_i^{\ast}=F_i,\ (q^{h})^{\ast}=q^{-h}\]
\[(resp.\ \omega(E_i)=F_i,\ \omega(F_i)=E_i,\ \omega(q^{h})=q^{-h}).\]
One can see that $\ast$ and $\omega$ can be extended to involutions
on $\widetilde{U}$, denoted by the same symbols, with
$(a_{\lambda})^{\ast}=\omega(a_{\lambda})=a_{-\lambda}$.

\subsection{Crystal basis and canonical basis of $\widetilde{U}$}

For $\lambda, \mu\in P^+$, let $V(\lambda)\otimes V(-\mu)$ be the
tensor product of irreducible integrable highest weight $U$-module
$V(\lambda)$ of highest weight $\lambda$ with irreducible integrable
lowest weight $U$-module $V(-\mu)$ of lowest weight $-\mu$. Note
that it is, by the coproduct of $U$, also a $U$-module and we denote
it also by $V(\lambda,-\mu)$. Let $u_{\lambda}$ (resp. $u_{-\mu}$)
be the highest (resp. lowest) weight vector of $V(\lambda)$ (resp.
$V(-\mu)$) and set $u_{\lambda,-\mu}=u_{\lambda}\otimes u_{-\mu}\in
V(\lambda,-\mu)$. It is known in \cite{Kashiwara:1994, Lusztig:1992} that $V(\lambda,-\mu)$ is a cyclic
$U$-module generated by $u_{\lambda,-\mu}$ and that it admits a crystal
basis
\[B(\lambda,-\mu)\triangleq B(\lambda)\otimes B(-\mu)\] where $B(\lambda)$ and $B(-\mu)$
are highest and lowest weight crystals respectively. The
corresponding global basis of $V(\lambda,-\mu)$ is constructed in
\cite{Lusztig:1992} and following Lusztig, we call it canonical basis, which
is denoted by $\{G(b)\ |\ b\in B(\lambda,-\mu)\}$.
\begin{definition}
\begin{itemize}

  \item[(i)]
  For a $U$-module $M$ with a canonical basis,
a subspace $N$  of $M$ is called nice or compatible with the canonical basis
if $N$ is spanned over $k$ by a part of the canonical basis of $M$.
  \item[(ii)]
 For $U$-modules $M$ and $N$ with canonical bases, a homomorphism of $U$-modules
$\phi:\ M\longrightarrow N$ is called nice or compatible with the
canonical bases if it maps the canonical base element of $M$ to
that of $N$ or to zero and if the images of  two distinct canonical base elements are distinct when they are both nonzero.
 \item[(iii)]
  For a $U$-module $M$ with a canonical basis,
a filtration or composition series of $M$ is called nice or compatible
with the canonical basis if any submodules in the filtration or composition series is nice.
\end{itemize}
\end{definition}

In \cite{Lusztig:1992}, the following stability property plays a key role in
the construction of the canonical basis of $\widetilde{U}$.

\begin{proposition} (\cite{Lusztig:1992})
For $\lambda, \mu, \theta\in P^+$, the map
$\pi_{\lambda,\mu,\theta}: V(\lambda+\theta,-\theta-\mu)\longrightarrow V(\lambda,-\mu)$
 which takes $xu_{\lambda+\theta,
-\theta-\mu}$ to $xu_{\lambda,-\mu}$ for all $x\in U$ is a nice surjective $U$-map.
\end{proposition}

We see from the proposition that there is an embedding of crystals
$B(\lambda,-\mu)\hookrightarrow
B(\lambda+\theta,-\theta-\mu)$ and note that it is strict \cite{Kashiwara:1994}. For
$\lambda,\mu\in P^+$, let $\Phi:
Ua_{\lambda-\mu}\longrightarrow V(\lambda,-\mu)$ be the $U$-map
 taking $a_{\lambda-\mu}$ to $u_{\lambda,-\mu}$.
It is known that $\widetilde{U}$ as well as each $Ua_{\lambda}$
have canonical bases and $\Phi$ is a nice surjective $U$-map. We
denote the crystal basis of $\widetilde{U}$ (resp. $Ua_{\lambda}$)
 by $\widetilde{B}$ (resp. $B(Ua_{\lambda})$). Hence we have the embedding of crystals
\[B(\lambda,-\mu)\hookrightarrow
 B(Ua_{\lambda-\mu}).\] It can be viewed as
$B(\lambda,-\mu)\subseteq B(\lambda+\theta,-\theta-\mu)\subseteq B(Ua_{\lambda-\mu})
 \subseteq \widetilde{B}.$ Note that $B(Ua_{\lambda})$ can be written as
 $B(\infty)\otimes T_{\lambda}\otimes B(-\infty)$ where $B(\pm\infty)$
 is the crystal basis of $U^{\mp}$ and $T_{\lambda}$ is a crystal consisting
  of a single element $t_{\lambda}$ with $\epsilon_i(t_{\lambda})=\phi_i(t_{\lambda})=-\infty$
 for all $i\in I$.
 For $b\in B(\lambda,-\mu)\subseteq\widetilde{B}$, we denote by the same $G(b)$ the
 corresponding canonical base element in $V(\lambda,-\mu)$ or $\widetilde{U}$ if this causes no confusion.
It is known that  $\ast$ induces a bijection on
$\widetilde{B}$ such that
$(b_1\otimes t_{\lambda}\otimes b_2)^{\ast}=b_1^{\ast}\otimes t_{-\lambda-wt(b_1)-wt(b_2)}\otimes b_2^{\ast}$ and
 $G(b)^{\ast}=G(b^{\ast})$ \cite{Kashiwara:1994}.

For any $\lambda\in P$, Kashiwara defined in \cite{Kashiwara:1994} an extremal weight $U$-module $V^{max}(\lambda)$
  which admits a crystal basis $B^{max}(\lambda)$ consisting of all $\ast$-extremal vectors in
   $B(Ua_{\lambda})$. We have
  $V^{max}(\lambda)\cong V^{max}(w\lambda)$ for any $w\in W$ and $V^{max}(\lambda)\cong V(\lambda)$ if
   $\lambda\in \pm P^+$. It is also known that for any connected
   component $B$ of $\widetilde{B}$,
   there is an $l>0$ such that $(wt(b),wt(b))\leqslant l$ for all $b\in B$. Moreover,
$B$ contains an extremal vector and can be embedded into $B^{max}(\mu)$ for some $\mu\in P$  \cite{Kashiwara:1994}.

For $\mathfrak g$ of affine type, let $c\in\mathfrak h$ be the canonical central element of $\mathfrak g$. Given
$\lambda\in P$, we define the level of $\lambda$ to be the integer $\lambda(c)$, denoted by $level(\lambda)$. Since
an integral weight $\lambda$ of positive (resp. negative) level is $W$-conjugate to a dominant (resp. anti-dominant) weight,
it follows from the previous paragraph that $B(Ua_{\lambda})$ is a union of highest (resp. lowest) weight crystals.

Denote by  $S^{\#}$  the cardinality of the set $S$. Given two
crystals $B_1$ and $B_2$ with $B_1$ connected, denote by $[B_2:
B_1]$  the cardinality of the set which consists of all connected
components of $B_2$ isomorphic to $B_1$, i.e.
\[[B_2: B_1]=\{B\subset B_2\ |\ B\cong B_1\}^{\#}.\]
The following result was proved in \cite{LZ:2009}.
\begin{proposition}\label{multiplicity}
  For $\lambda\in P^+$ and $\mu\in P$, $[B(Ua_{\mu}):
B(\lambda)]=dimV(\lambda)_{\mu}.$
\end{proposition}

\subsection{Bilinear Form on $\widetilde{U}$}

We introduce another anti-automorphism of $U$, denoted
by $\Psi$ \cite{Kashiwara:1994}, such that
\[\Psi(E_i)=q_i^{-1}t_i^{-1}F_i,\ \Psi(F_i)=q_i^{-1}t_iE_i,\ \Psi(q^h)=q^h,\ \Psi(q)=q,\]
where $q_i=q^{\frac{(\alpha_i,\alpha_i)}{2}}$ and $t_i=q^{\frac{(\alpha_i,\alpha_i)}{2}h_i}$. One can easily check that
$\Psi^2=id$ and $\Psi$ commutes with $\omega$ introduced previously.

For $\lambda\in \pm P^+$, there is a unique non-degenerate symmetric bilinear form $(\ ,\ )$ on $V(\lambda)$ such that
\[(u_{\lambda}, u_{\lambda})=1\ \textrm{and} \ (Pu,v)=(u,\Psi(P)v)\ \textrm{for\ all}\ P\in U,\ u,v\in V(\lambda).\]
See that
\[(Pu_{\lambda}, Qu_{\lambda})=(\omega(P)u_{-\lambda}, \omega(Q)u_{-\lambda})\ \ \ \ \ \ \ \ \ \ \ \ (2.1).\]
Indeed, one can define a bilinear form $((\ ,\ ))$ on $V(\lambda)$ by
\[((Pu_{\lambda},Qu_{\lambda}))=(\omega(P)u_{-\lambda},\omega(Q)u_{-\lambda}),\] and (2.1) follows from the
uniqueness of the bilinear form on $V(\lambda)$.
Given $\lambda, \mu\in P_+$, we define a symmetric bilinear form $(\ ,\ )$ on $V(\lambda,-\mu)$ by
$(u_1\otimes v_1, u_2\otimes v_2)=(u_1,v_1)(u_2,v_2).$ Since $\Psi$ commutes with the coproduct $\Delta$, i.e.
$(\Psi\otimes\Psi)\Delta=\Delta\Psi$, it implies that
$(Pu,v)=(u,\Psi(P)v)$ for all $P\in U,\ u,v\in V(\lambda,-\mu)$.

\begin{lemma} (\cite{Kashiwara:1994})
For $P, Q\in U$ and $\theta\in P$, there exists
a unique polynomial $f(x)$ in $x=(x_i)_i\in I$ such that for any $\lambda, \mu\in P^+$ with $\lambda-\mu=\theta$,
$(Pu_{\lambda,-\mu}, Qu_{\lambda,-\mu})=f(x)\ \textrm{with}\ x_i=q_i^{\lambda(h_i)}.$
\end{lemma}

The bilinear form on $Ua_{\theta}$ is then defined by $(Pa_{\theta}, Qa_{\theta})=f(0)$ and this
extends to a bilinear form on $\widetilde{U}$ such
  that $(Ua_{\theta_1}, Ua_{\theta_2})=0\ \textrm{for}\  \theta_1\neq \theta_2.$ It was shown in
  \cite{Kashiwara:1994} that $(\ ,\ )$ on $\widetilde{U}$ is symmetric and it satisfies
  \[(u,v)=(u^{\ast},v^{\ast})\ \textrm{and}\ (Pu,v)=(u,\Psi(P)v)\ \textrm{for\ all}\ P\in U, u,v\in \widetilde{U}.\]
Let $A_0$ be the subring of $k$ consisting of all rational functions regular at $q=0$ and $\widetilde{U}_{\Bbb Z}$ be the
$\Bbb Z$-form of $\widetilde{U}$ \cite{Kashiwara:1994}. The crystal lattice $L(\widetilde{U})$ over $A_0$ and the canonical basis of $\widetilde{U}$ are
characterized by the bilinear form
\cite{Kashiwara:1994}.
\begin{proposition} \label{bilinearform}
\begin{itemize}
  \item[(i)]
   $L(\widetilde{U})=\{x\in\widetilde{U}\ |\ (u,u)\in A_0\}$
    \item[(ii)]
   If $u\in\widetilde{U}_{\Bbb Z}$ and $(u,u)\in 1+qA_0$, then
$u\equiv G(b)$ or $-G(b)$ mod $qL(\widetilde{U})$ for some $b\in\widetilde{U}$.
\end{itemize}
\end{proposition}

  We define another bilinear form $((\ ,\ ))$ on $\widetilde{U}$ by
  \[((u,v))=(\omega(u),\omega(v))\ \textrm{for\ all}\ u,v\in\widetilde{U}.\]

\begin{proposition} \label{invariant}
$((u,v))=(u,v)$ for all $u,v\in\widetilde{U}$.
\end{proposition}
\begin{proof}
Assume that $u=Pa_{\theta}$, $v=Qa_{\theta}$ where $P,Q\in U$, $\theta\in P$. Let
$\Psi(P)Q=\sum_{j}x_j^+x_j^-$ where $x_j^{\pm}\in U^{\pm}\otimes k[q^{h}:h\in P^{\vee}]$. We have
\begin{equation*}
\begin{split}
((Pa_{\theta},Qa_{\theta})) &= (\omega(P)a_{-\theta},\omega(Q)a_{-\theta}) =(a_{-\theta},\Psi(\omega(P))\omega(Q)a_{-\theta})\\
                            &= (a_{-\theta},\omega(\Psi(P)Q)a_{-\theta})=\sum_{j}(a_{-\theta},\omega(x_j^+x_j^-)a_{-\theta}) \\
                            &= \sum_{j}(\omega(\Psi(x_j^+))a_{-\theta},\omega(x_j^-)a_{-\theta}) =\sum_j((\Psi(x_j^+)a_{\theta},x_j^-a_{\theta})) \\
\end{split}
\end{equation*}
Since $q^{h}a_{\theta}=q^{\theta(h)}a_{\theta}$, it is sufficient to show the equality when $P, Q\in U^-$. Given $\lambda, \mu\in P^+$ such
that $\lambda-\mu=\theta$,
\[(Pu_{\lambda,-\mu}, Qu_{\lambda,-\mu})=(Pu_{\lambda}\otimes u_{-\mu},Qu_{\lambda}\otimes u_{-\mu})=(Pu_{\lambda},Qu_{\lambda})=f(x)\]
Meanwhile we have,
\[(\omega(P)u_{\mu,-\lambda}, \omega(Q)u_{\mu,-\lambda}))=
(u_{\mu}\otimes
\omega(P)u_{-\lambda},u_{\mu}\otimes\omega(Q)u_{-\lambda})\]
\[=(\omega(P)u_{-\lambda},\omega(Q)u_{-\lambda})=(Pu_{\lambda},Qu_{\lambda})=f(x).\]
Hence $((Pa_{\theta},Qa_{\theta}))=f(0)=(Pa_{\theta},Qa_{\theta}).$
\end{proof}

\begin{corollary}\label{omega}
\begin{itemize}
  \item[(i)]
   $\omega(L(\widetilde{U}))=L(\widetilde{U})$.
  \item[(ii)]
   $\omega(\widetilde{B})=\widetilde{B}$.
  \item[(iii)]
   $\omega(G(b))=G(\omega(b))$ for all $b\in\widetilde{U}$.
\end{itemize}
\end{corollary}
\begin{proof}
 $(\romannumeral 1)$ follows immediately from Proposition \ref{bilinearform} and \ref{invariant}. The proof
 of  $(\romannumeral 3)$ is similar to \cite{Kashiwara:1991}. We only show  $(\romannumeral 2)$ here. Given
 $b=b_1\otimes t_{\lambda}\otimes b_2\in\widetilde{B}$ with $ht(wt(b_i))=l_i$, i=1, 2, we have
 \[G(b)\equiv G(b_1)G(b_2)a_{\lambda}\ \textrm{mod}\ F_{l_1-1}(U^-)F_{l_2-1}(U^+)a_{\lambda}\ \ \ \ \ \ \ \ \ \ \ \ (2.2).\]
Since $\omega: U^{\pm}\longrightarrow U^{\mp}$ induces $\omega : B(\mp\infty)\longrightarrow B(\pm\infty)$ and
it maps canonical base elements of $U^{\pm}$ to those of $U^{\mp}$, applying $\omega$ to (2.2) we have
\[\omega(G(b))\equiv \omega(G(b_1))\omega(G(b_2))a_{-\lambda}=G(\omega(b_1))G(\omega(b_2))a_{-\lambda}\]
\[\equiv G(\omega(b_2))G(\omega(b_1))a_{-\lambda}\ \textrm{mod}\ F_{l_2-1}(U^-)F_{l_1-1}(U^+)a_{-\lambda}.\]
Since $\omega(G(b))=G(b^{\prime})$ or $-G(b^{\prime})$ for some $b^{\prime}\in\widetilde{B}$, we obtain that
$b^{\prime}=\omega(b_2)\otimes t_{-\lambda}\otimes \omega(b_1)$ and $\omega(G(b))=G(b^{\prime})$.
\end{proof}

\section{A quotient algebra of $\widetilde{U}$}

 Throughout this
section, a pair of dominant weights $(\lambda, \mu)$ is fixed.

\subsection{Filtration}
In this subsection, we  recall the construction of some nice filtrations of the
$U$-module $V(\lambda,-\mu)$ in \cite{LZ:2009}. In order to
obtain nice submodules of $V(\lambda, -\mu)$, we need the following lemma due to Kashiwara \cite{Kashiwara:1993}. See
also \cite{LZ:2009} for more details.
\begin{lemma}\label{kashiwara}
\begin{itemize}
  \item[(i)]
Let $M$ be an integrable
$U$-module
with a canonical basis. If $N$ is a nice $U^+$-submodule
of $M$, then $UN=U^-N$ is a nice $U$-submodule of $M$. More
precisely,  $UN=\bigoplus_{b\in B(UN)\subseteq B(M)} kG(b).$
  \item[(ii)]
  $B(UN)=\{\tilde{f}_{i_1}\cdots\tilde{f}_{i_m}b\ |\  m\geqslant 0,
i_1,\cdots,i_m\in I,b\in B(N) \}\setminus \{0\}$.
\end{itemize}
\end{lemma}

One can define a total order $<$ on the lowest weight crystal $B(-\mu)$ such that $b_1<b_2$
if $wt(b_1)<wt(b_2)$. Indeed, for $b\in B(-\mu)$, set
$l(b)=m$ if $b$ is of the form $\tilde{e}_{i_1}\cdots\tilde{e}_{i_m}u_{-\mu}$. It is well-defined since
$b=\tilde{e}_{i_1}\cdots\tilde{e}_{i_m}u_{-\mu}=\tilde{e}_{j_1}\cdots\tilde{e}_{j_l}u_{-\mu}$
implies $l=m$ by comparing the weights. We arrange the order on  $I$
as
\[1<2<\cdots<n.\] Let $|b|$ denote the $l(b)$-tuple $(i_1, \cdots,
i_{l(b)})$ such that $(i_1, \cdots, i_{l(b)})$ is minimal in
lexicographic order among tuples $(j_1, \cdots,
 j_{l(b)})$ such that
 $\tilde{e}_{j_1}\cdots\tilde{e}_{j_{l(b)}}u_{-\mu}=b$,
 i.e. $|b|=min\{ (j_1, \cdots, j_{l(b)})\ |\ b=\tilde{e}_{j_1}
 \cdots\tilde{e}_{j_{l(b)}}u_{-\mu} \}.$
 Set $|u_{-\mu}|=0$. The order on $B(-\mu)$ is then defined as
 follows,
 \[b_1\leqslant b_2\ \textrm{iff} \ l(b_1)<l(b_2)\ \textrm{or}\  l(b_1)=l(b_2)\
\textrm{but} \ |b_1|\leqslant|b_2|. \]

For $b\in B(-\mu)$, we define a subspace $V_b(-\mu)$ of $V(-\mu)$ as
 \[V_b(-\mu)\triangleq \sum_{c\geqslant b}kG(c)\] which is easily shown to
 be a $U^+$-submodule. Hence $u_{\lambda}\otimes V_b(-\mu)$ is a
 $U^+$-submodule of $V(\lambda,-\mu)$ which
  has a basis $\{u_{\lambda}\otimes G(c)\ |\ c\geqslant b,\ c\in B(-\mu) \}.$
  Since $u_{\lambda}\otimes G(c)=G(u_{\lambda}\otimes c)$, $u_{\lambda}\otimes V_b(-\mu)$ is actually a
 nice $U^+$-submodule. Define
 $F_{\lambda}(b)$ to be a $U$-submodule of $V(\lambda, -\mu)$ generated by $u_{\lambda}\otimes V_b(-\mu)$, i.e.
 \[F_{\lambda}(b)=U(u_{\lambda}\otimes V_b(-\mu)).\]
It follows from Lemma~\ref{kashiwara} that
$F_{\lambda}(b)$ is a nice $U$-submodule of
$V(\lambda,-\mu)$ and
\[B(F_{\lambda}(b))=\{\tilde{f}_{i_1}\cdots\tilde{f}_{i_m}(u_{\lambda}\otimes
c)\ |\
i_1,\cdots,i_m\in I,\ c\in B(-\mu),\ c\geqslant b\}\setminus\{0\}.\]
Moreover, by comparing the crystal basis, we have in \cite{LZ:2009} the following result.

\begin{theorem}
For two neighbors $b<c\in B(-\mu)$,
$F_{\lambda}(b)/F_{\lambda}(c)\cong V(\lambda+wt(b))$ if
$\tilde{e}_i(u_{\lambda}\otimes b)=0$ for all $i\in I$, otherwise
$F_{\lambda}(b)=F_{\lambda}(c)$.
\end{theorem}

Hence we get a nice descending filtration of $V(\lambda, -\mu)$
\[V(\lambda, -\mu)=F_{\lambda}(b_1)\supseteq F_{\lambda}(b_2)\supseteq F_{\lambda}(b_3)
\supseteq\cdots\ \ \ \ \ \ \ \ \ \ \ \ \ \ \ (3.1)\] where $u_{-\mu}=b_1<b_2<b_3<\cdots$
is a complete list of $B(-\mu)$.

\begin{remark}
\begin{itemize}
    \item[(i)]
    In our construction, the total order on
$B(-\mu)$ is fixed while in fact we can choose any total
  order such that $b_1<b_2$ if  $wt(b_1)<wt(b_2)$.
   \item[(ii)]
   Similarly one can also define a total order
on $B(\lambda)$ such that $b_1< b_2$ if $wt(b_1)<wt(b_2)$. Set
\[F^{-\mu}(b)\triangleq\sum_{c\leqslant b}U(G(c)\otimes u_{-\mu})\]
with which we can also construct a nice filtration of $V(\lambda,-\mu)$
where the quotient of two neighbors is isomorphic either to
an irreducible lowest weight module or to 0.
\end{itemize}
\end{remark}

Let $W(\lambda,-\mu)$ be a subspace of $V(\lambda,-\mu)$ defined by
\[W(\lambda,-\mu)\triangleq\bigcap_{b\in B(-\mu)}F_{\lambda}(b).\]
Set $M(\lambda,-\mu)=V(\lambda,-\mu)/W(\lambda,-\mu)$. Denote by
$B^{\prime}$ (resp. $B^{\prime}(\lambda,-\mu)$) the sub-crystal of
$\widetilde{B}$ (resp. $B(\lambda, -\mu)$) which is a union of all
connected components of $\widetilde{B}$ (resp. $B(\lambda, -\mu)$)
that are not highest weight crystals. We have the following proposition in \cite{LZ:2009}.
\begin{proposition}
\begin{itemize}
\item[(i)]
$W(\lambda,-\mu)$ is a nice $U$-submodule of
$V(\lambda,-\mu)$ and $B(W(\lambda,-\mu))=
B^{\prime}(\lambda,-\mu)$.
\item[(ii)]
$M(\lambda,-\mu)$ admits a canonical basis and
$B(M(\lambda,-\mu))=
B(\lambda,-\mu)\setminus B^{\prime}(\lambda,-\mu).$
\end{itemize}
\end{proposition}

\begin{remark}
One can see that $U(\lambda,-\mu)\triangleq\bigcap_{b\in B(\lambda)}F_{-\mu}(b)$ has a crystal basis $B^{\prime\prime}(\lambda,\mu)$ as well as
a canonical basis where $B^{\prime\prime}(\lambda,-\mu)$ consists of all connected components of $B(\lambda,-\mu)$ that are not lowest weight crystals.
Similarly $N(\lambda,-\mu)\triangleq V(\lambda,-\mu)/U(\lambda,-\mu)$ admits a canonical basis.
\end{remark}

 Note that when $\mathfrak g$ is of finite type, $V(\lambda,-\mu)$ is finite
dimensional. Hence there are finitely many terms in the filtration
(3.1) and furthermore, we can obtain a nice composition series of
$V(\lambda,-\mu)$ \cite{LZ:2009} by deleting the superfluous terms in
(3.1) which provides a complete proof to the conjecture raised by
Lusztig \cite{Lusztig:1992}. Moreover, $W(\lambda, -\mu)=0$ and
$M(\lambda,-\mu)=V(\lambda,-\mu)$ in this case. But when $\mathfrak
g$ is of affine or indefinite type, the situation is quite
different. For $\mathfrak g$ of affine type, the following result
was shown in \cite{LZ:2009}.
\begin{proposition} \label{affine}
\begin{itemize}
\item[(i)]
$W(\lambda,-\mu)=N(\lambda,-\mu)=0$ and
$M(\lambda,-\mu)=U(\lambda,-\mu)$
$=V(\lambda, -\mu)$ if
$level(\lambda-\mu)>0$.
\item[(ii)]
$W(\lambda,-\mu)=N(\lambda,-\mu)=V(\lambda, -\mu)$ and
$M(\lambda,-\mu)=U(\lambda,-\mu)=0$ if $level(\lambda-\mu)<0$.
\item[(iii)]
$M(\lambda,-\mu)=N(\lambda,-\mu)$ is a 1-dimensional trivial module
if $\lambda-\mu\in P_0$, otherwise if $\lambda-\mu\notin P_0$ is of
level 0, $W(\lambda,-\mu)=U(\lambda,-\mu)=V(\lambda, -\mu)$
and $M(\lambda,-\mu)=N(\lambda,-\mu)=0$.
\end{itemize}
\end{proposition}

\subsection{$\mathfrak{U}$}

We denote by $\mathcal O^+$ (resp. $\mathcal O^-$) the completely reducible category whose objects are
direct sums of irreducible integrable highest (resp. lowest) weight $U$-modules. Note that $\mathcal O^+$ here is often referred to as
$\mathcal O_{int}$ in other literatures.

\begin{theorem} \label{equvilant} For $b\in \widetilde{B}$, the
following conditions are
equivalent.
\begin{itemize}
\item[(i)]
$G(b)$ acts on $V(\lambda)$ as zero for all
$\lambda\in P^+$.
\item[(ii)]
$G(b)$ acts on $M$ as zero for any
$M\in ob(\mathcal O^+)$.
\item[(iii)]
$b\in B^{\prime}$.
\end{itemize}
\end{theorem}
\begin{proof}
The equivalence of $(\romannumeral 1)$ and
$(\romannumeral 2)$ is clear. If $b$ satisfies $(\romannumeral 2)$,
we show that it satisfies $(\romannumeral 3)$. Otherwise assume that
$b\notin B^{\prime}$, $b$ is contained in a highest weight
subcrystal of $\widetilde{B}$. There exist $\lambda$, $\mu\in P^+$
such that
$b\in B(\lambda, -\mu)\subset \widetilde{B}.$
We rewrite the nice filtration (3.1) of $V(\lambda, -\mu)$ as
\[V(\lambda, -\mu)=F_0\supseteq F_1\supseteq\cdots\supseteq
F_l\supseteq\cdots\ \ \ \ \ \ \ \ \ \ \ \ \ \ \ (3.2).\] There exists an $s\geqslant
0$ such that $G(b)\in F_s$ but $G(b)\notin F_{s+1}$. Hence
\[0\neq G(b)(u_{\lambda,-\mu}+F_{s+1})\in V(\lambda,
-\mu)/F_{s+1}\] where $V(\lambda, -\mu)/F_{s+1}$ is an object in
$\mathcal O^+$. This contradicts $(\romannumeral 2)$. Finally we
show that $(\romannumeral 3)$ implies $(\romannumeral 1)$. Assume
that $G(b)V(\lambda^{\prime})\neq 0$ for some $\lambda^{\prime}\in
P^+$, then there exists an $m\in V(\lambda^{\prime})_{\xi}$ such
that $G(b)m\neq 0$ and $b\in B(Ua_{\xi})\subset \widetilde{B}$. We
can find $\lambda, \mu\in P^+$ with $\lambda-\mu=\xi$ such that
$b\in B(\lambda, -\mu)\subset B(Ua_{\xi})$ and there exists a
homomorphism of $U$-modules $\phi: V(\lambda, -\mu)\longrightarrow
V(\lambda^{\prime})$ which takes $u_{\lambda,-\mu}$ to $m$. Since
$b\in B^{\prime}$, $G(b)\in W(\lambda, -\mu)$, that is, $G(b)\in
F_s$ for any $F_s$ in the filtration (3.2). Restricting $\phi$ on
$F_s$, we get a $U$-morphism
\[\phi|_{F_s}: F_s\longrightarrow V(\lambda^{\prime}).\]
Since the set of generators of $F_s$ is of the form
$u_{\lambda}\otimes V_b(-\mu)$ for some $b$, the corresponding
weights of these generators are not lower than or equal to $\lambda^{\prime}$
for a sufficient large $s$ by the construction. Hence $\phi|_{F_s}$
is zero for $s>>0$. It follows that $\phi(G(b))=G(b)m=0$ which is a
contradiction.
\end{proof}

 As is known in \cite{Joseph:1995}, if $u\in U$ acts on each $M\in ob(\mathcal O^+)$ as
 zero, then $u=0$. But it is not true for $u\in \widetilde{U}$ and
 $\mathfrak g$ of affine or indefinite type by the above
 theorem.
 \begin{proposition} \label{annihilate}
 For $u=\sum k_bG(b)\in \widetilde{U}$ such that $u$ acts on
 $M$ as zero for all $M\in ob(\mathcal O^+)$ and if $k_b\neq 0$, then $b\in
 B^{\prime}$.
 \end{proposition}

\begin{proof}
We assume that $k_{b_0}\neq 0$ for some $b_0\notin
 B^{\prime}$. There exist $\lambda, \mu\in P^+$ such that
$b_0\in B(\lambda, -\mu)\subset \widetilde{B}.$ Since $b_0\notin
 B^{\prime}$, there exists an $s$ such that $G(b_0)\in F_s$ but $G(b_0)\notin
 F_{s+1}$ where $F_s$ and $F_{s+1}$ are in the filtration (3.2) of
 $V(\lambda, -\mu)$. Hence we have
 \[0\neq u(u_{\lambda,-\mu}+F_{s+1})\in V(\lambda,
 -\mu)/F_{s+1}\] with $V(\lambda,
 -\mu)/F_{s+1}\in ob(\mathcal O^+)$. This is a contradiction.
\end{proof}
 By this proposition we know that any $u\in\widetilde{U}$ annihilating all $M\in ob(\mathcal O^+)$ is a linear combination of $G(b)$'s
 with $b\in B^{\prime}$. Denote by $\widetilde{U}^{\prime}$ the set of
 all such $u$'s. It follows from Theorem~\ref{equvilant} and Proposition~\ref{annihilate} that
 \begin{theorem}
 $\widetilde{U}^{\prime}$ is a nice two-sided ideal of
 $\widetilde{U}$ and it admits a crystal basis $B^{\prime}$.
 \end{theorem}

 We define  $\mathfrak{U}$ to be the quotient of $\widetilde{U}$ by
 $\widetilde{U}^{\prime}$, i.e.
 $\mathfrak{U}\triangleq\widetilde{U}/\widetilde{U}^{\prime}.$ Hence $\mathfrak{U}$ inherits from $\widetilde{U}$
 a canonical basis and we denote by $\mathfrak{B}$ the corresponding crystal basis. One can see from the definition
 of $B^{\prime}$ and $\mathfrak{B}=\widetilde{B}\setminus B^{\prime}$ that $\mathfrak{B}$ is a union of all highest weight sub-crystals of $\widetilde{B}$.
 We know also from Theorem~\ref{equvilant} that any $M\in ob(\mathcal O^+)$ is also a representation
of $\mathfrak{U}$. Note that when $\mathfrak g$ is of finite type,
$\mathfrak{U}=\widetilde{U}$. If $\mathfrak g$ is of affine type, it follows from
Proposition~\ref{affine} that $\mathfrak{U}$ is isomorphic to the subalgebra of $\widetilde{U}$
generated by $Ua_{\xi}$ and $a_{\eta}$ for all $\xi$ with a
positive level and $\eta\in P_0$.

\begin{remark}
Similarly one can define $\widetilde{U}^{\prime\prime}$ to be the set of all $u\in\widetilde{U}$ such that
$u$ annihilates all $M\in ob(\mathcal O^-)$. Then $\widetilde{U}^{\prime\prime}$ is also a nice ideal of $\widetilde{U}$ with
a crystal basis $B^{\prime\prime}$ where $B^{\prime\prime}$ consists of all connected components of $\widetilde{B}$ that are
not lowest weight crystals. We denote by $\mathfrak{V}$ the quotient algebra $\widetilde{U}/\widetilde{U}^{\prime\prime}$ which admits
both a cystal basis and a canonical basis.
\end{remark}

\section{Cells in $\mathfrak{U}$}

Recall that we define $M(\lambda,-\mu)$ which is a representation of
$U$ as well as $\widetilde{U}$ in the previous section. Also it can be
viewed as a representation of $\mathfrak{U}$. To see that, we need
the following lemma.
\begin{lemma}
 For $b\in B^{\prime}\subset B$ and
$\lambda, \mu\in P^+$, $G(b)M(\lambda, -\mu)=0$.
\end{lemma}
\begin{proof}
 We only show that $G(b)V(\lambda,
-\mu)\subseteq W(\lambda, -\mu)$. Since $V(\lambda, -\mu)/F_s\in
ob(\mathcal O^+)$ for any $F_s$ in the filtration (3.2), by
Theorem~\ref{equvilant} we have \[G(b)(V(\lambda, -\mu)/F_s)=0\] which
means $G(b)V(\lambda, -\mu)\subseteq F_s$. Hence we have
\[G(b))V(\lambda, -\mu)\subseteq \bigcap_{s\geqslant 0}F_s=W(\lambda, -\mu).\]
\end{proof}
Applying this lemma we have $\widetilde{U}^{\prime}M(\lambda, -\mu)=0$ and thus we
equip $M(\lambda, -\mu)$ with a $\mathfrak{U}$-action. Quotient by
$W(\lambda, -\mu)$, we obtain from (3.2) a filtration of $M(\lambda,
-\mu)$ consisting of nice $U$ or $\mathfrak{U}$-submodules
\[M(\lambda, -\mu)=M_0\supseteq M_1\supseteq\cdots\supseteq
M_l\supseteq\cdots\ \ \ \ \ \ \ \ \ \ \ \ \ (4.1)\] where
$M_i=F_i/W(\lambda, -\mu).$ Denote by $v_{\lambda, -\mu}$ the image
of $u_{\lambda,-\mu}$ in $M(\lambda, -\mu)$. Hence the map
\[\bar{\alpha}_{\lambda,-\mu}:\ \mathfrak{U}\longrightarrow M(\lambda, -\mu)\ \ \
x\longmapsto xv_{\lambda, -\mu}\] takes the canonical base elements
of $\mathfrak{U}$ to those of $M(\lambda, -\mu)$ or to zero. For $\xi\in
P^+$, let $M(\lambda, -\mu)_{[\xi]}$ ( resp. $\mathfrak{U}_{[\xi]}$) be
the subspace of $M(\lambda, -\mu)$ ( resp. $\mathfrak{U}$) spanned by
all $G(b)$'s such that $b$ is contained in a subcrystal of
$B(M(\lambda, -\mu))$ ( resp. $\mathfrak{B}$) isomorphic to
$B(\xi)$. Note that $M(\lambda, -\mu)_{[\xi]}$ is usually not a
$U$-submodule of $M(\lambda, -\mu)$. Set $M(\lambda,
-\mu)_{[\geqslant\xi]}$, $M(\lambda, -\mu)_{[>\xi]}$,
$\mathfrak{U}_{[\geqslant\xi]}$ and $\mathfrak{U}_{[>\xi]}$ as follows
\[M(\lambda, -\mu)_{[\geqslant\xi]}\triangleq\bigoplus_{\eta\geqslant\xi}M(\lambda,
-\mu)_{[\eta]},\ \
\mathfrak{U}_{[\geqslant\xi]}\triangleq\bigoplus_{\eta\geqslant\xi}\mathfrak{U}_{[\eta]},\]
\[M(\lambda, -\mu)_{[>\xi]}\triangleq\bigoplus_{\eta>\xi}M(\lambda, -\mu)_{[\eta]},\ \
\mathfrak{U}_{[>\xi]}\triangleq\bigoplus_{\eta>\xi}\mathfrak{U}_{[\eta]}.\]

For a $U$-module $M\in ob(\mathcal O^+)$ with a canonical basis, $M$ can
be written as $M=\bigoplus_{\lambda\in P^+}M[\lambda]$ where
$M[\lambda]$ is the sum of all submodules of $M$ isomorphic to
$V(\lambda)$. Here $M[\lambda]$ is usually not nice. But it is known
in \cite{Lusztig:1993} that $M[\xi]$ is a nice $U$-submodule of $M$ for a maximal $\xi$,
i.e. $\xi$ is maximal in the sense of dominant order among all
$\lambda$ such that $M[\lambda]\neq 0$. Moreover, both
$M[\geqslant\xi]\triangleq\bigoplus_{\lambda\geqslant\xi}M[\lambda]$
and $M[>\xi]\triangleq\bigoplus_{\lambda>\xi}M[\lambda]$ are nice.
In particular, when $\mathfrak g$ is of finite type,
$M(\lambda, -\mu)_{[\geqslant\xi]}=M(\lambda,
-\mu)[\geqslant\xi]$.\\

Similar to \cite{Lusztig:1993} we have the following lemma.
\begin{lemma}\label{equvilant2}
For $x\in\mathfrak{U}$
and $\xi\in P^+$, the following are
equivalent
\begin{itemize}
\item[(i)]
  $x\in\mathfrak{U}_{[\geqslant\xi]}$.
\item[(ii)]
For all $\lambda, \mu\in P^+$,
$xv_{\lambda, -\mu}\in M(\lambda,-\mu)_{[\geqslant\xi]}$.
\item[(iii)]
For any $M\in ob(\mathcal O^+)$ and $m\in M$, $xm\in M[\geqslant\xi]$.
\item[(iv)]
If $x$ acts on $V(\eta)$ as a nonzero map for some $\eta\in P^+$, then
$\eta\geqslant\xi$.
\end{itemize}
\end{lemma}

\begin{proof}
It is clear that the equivalence of
$(\romannumeral 1)$ and $(\romannumeral 2)$ follows from
definitions of  $\mathfrak{U}_{[\geqslant\xi]}$ and $M(\lambda,-\mu)_{[\geqslant\xi]}$.
 $(\romannumeral 3)$ and $(\romannumeral 4)$ are
equivalent since any $M\in ob(\mathcal O^+)$ can be written as a direct sum
of some $V(\eta)$'s. If $x$ satisfies $(\romannumeral 3)$£¬ we show
that it satisfies $(\romannumeral 2)$. Set $x=\sum k_bG(b)$. Assume
that $(\romannumeral 2)$ does not hold, then there exists some
$b_0\in B(M(\lambda, -\mu))\subseteq \mathfrak{B}$ with $k_{b_0}\neq
0$ such that $b_0$ is contained in a subcrystal of $\mathfrak{B}$
isomorphic to $B(\eta)$ with $\eta\ngeqslant\xi$. It follows that
there exists an $s$ such that $G(b_0)\in M_s$ but $G(b_0)\notin
M_{s+1}$ where $M_s$ and $M_{s+1}$ are in the filtration (4.1). Thus
\[M_s/M_{s+1}\cong V(\eta).\] Set $M\triangleq
M(\lambda,-\mu)/M_{s+1}$ and $m=v_{\lambda,-\mu}+M_{s+1}\in M$. Then
$M\in ob(\mathcal O^+)$ and $xm\notin M[\geqslant\xi]$ which contradicts
$(\romannumeral 3)$. Conversely we show that $(\romannumeral 2)$
implies $(\romannumeral 3)$. For any $M\in ob(\mathcal O^+)$ and $m\in
M_{\theta}$, there exists $\lambda, \mu\in P^+$ with
$\lambda-\mu=\theta$ such that $xv_{\lambda,-\mu}\neq 0$ and $\phi:
V(\lambda, -\mu)\longrightarrow M,\ u_{\lambda,-\mu}\longmapsto m$
is a nonzero $U$-map. As in the proof of Theorem~\ref{equvilant}, one can see
that $\phi(W(\lambda, -\mu))=0$. Hence we have \[\bar{\phi}:\
M(\lambda, -\mu)\longrightarrow M,\ v_{\lambda, -\mu}\longmapsto m\]
a homomorphism of both $U$-modules and $\mathfrak{U}$-modules. As is proved before,
there exists an $s$ such that the weights
of the generators of $M_s$ are not lower than or equal to
any weight in $M$. Hence $\bar{\phi}(M_s)=0$ and furthermore,
$\bar{\phi}$ factors through the $\mathfrak{U}$-map
$\bar{\phi}^{\prime}:\ M(\lambda,-\mu)/M_s\longrightarrow M,\ \
v_{\lambda,-\mu}+M_s\longmapsto m.$ Since $xv_{\lambda, -\mu}\in
M(\lambda,-\mu)_{[\geqslant\xi]}$ and
$M(\lambda,-\mu)/M_s\in ob(\mathcal O^+)$, \[x(v_{\lambda,-\mu}+M_s)\in
(M(\lambda,-\mu)/M_s)[\geqslant\xi].\] It follows that
$\bar{\phi}^{\prime}(x(v_{\lambda,-\mu}+M_s))=xm\in M[\geqslant\xi]$
which proves $(\romannumeral 3)$.
\end{proof}

Similarly one can prove the following lemma since
\[\mathfrak{U}_{[>\xi]}=\sum_{\eta>\xi}\mathfrak{U}_{[\geqslant\eta]},\
M_{[>\xi]}=\sum_{\eta>\xi}M_{[\geqslant\eta]},\
M(\lambda,-\mu)_{[>\xi]}=\sum_{\eta>\xi}M(\lambda,-\mu)_{[\geqslant\eta]}.\]
\begin{lemma}\label{equvilant3}
For $x\in\mathfrak{U}$ and $\xi\in P^+$, the following are
equivalent
\begin{itemize}
\item[(i)]
$x\in\mathfrak{U}_{[>\xi]}$.
\item[(ii)]
For all $\lambda, \mu\in P^+$,
$xv_{\lambda, -\mu}\in M(\lambda,-\mu)_{[>\xi]}$.
\item[(iii)]
For any $M\in ob(\mathcal O^+)$ and $m\in M$, $xm\in M[>\xi]$.
\item[(iv)]
If $x$ acts on $V(\eta)$ as a nonzero map for
some $\eta\in P^+$, then $\eta>\xi$.
\end{itemize}
\end{lemma}

The corollary below follows immediately from Lemma~\ref{equvilant2} and Lemma~\ref{equvilant3}.
\begin{corollary}\label{ideal}
Both $\mathfrak{U}_{[\geqslant\xi]}$ and $\mathfrak{U}_{[>\xi]}$ are nice
two-sided ideals of $\mathfrak{U}$ for any $\xi\in P^+$.
\end{corollary}

\begin{remark}
For $\xi\in P^+$, we can define $\mathfrak{V}_{[\leqslant -\xi]}$ (resp.
$\mathfrak{V}_{[<-\xi]}$) to be the subset of $\mathfrak{V}$ consisting of
all $x$ such that $\eta\geqslant\xi$ (resp. $\eta>\xi$) if $x$ acts
on $V(-\eta)$ as a nonzero map. Similarly both of them are nice
ideals of $\mathfrak{V}$.
\end{remark}

For an integrable left $U$-module $M$ with finite dimensional weight
spaces, let $M^{\circ}$ denote the graded dual of $M$, i.e.
$M^{\circ}=\bigoplus_{\theta\in P} M_{\theta}^{\ast}\ \textrm{where}\ M=\bigoplus_{\theta\in P}M_{\theta}.$
Then there is a right $U$-action on $M^{\circ}$ as
\[(f\cdot x)(v)=f(xv)\ \textrm{for}\ f\in M^{\circ},\ v\in M, x\in U.\]
For instance, $V(\lambda)^{\circ}$ is an irreducible integrable
right $U$-module with highest weight $\lambda\in P^+$. Given a right
$U$-module $M$, we denote by ${}^{\ast}M$ the same $k$-vector space equipped with a left $U$-action as
\[x\circ m=m\cdot x^{\ast}\ \textrm{for}\ x\in U,\ m\in {}^{\ast}M.\]
It is clear that ${}^{\ast}V(\lambda)^{\circ}\cong V(-\lambda)$ as left $U$-modules for $\lambda\in \pm P^{+}$.
Given a left $U$-module $N$, define ${}^{\omega}N$ to be the left $U$-module with the underlying space
${}^{\omega}N=N$ such that
\[x\circ v=\omega(x)\cdot v\ \textrm{for}\ x\in U,\ v\in {}^{\omega}N.\] See that ${}^{\omega}V(\lambda)\cong V(-\lambda)$
for $\lambda\in \pm P^{+}$.
\begin{lemma}\label{star}
\begin{itemize}
\item[(i)]
Both $\ast$ and $\omega$ on $\widetilde{U}$ induce  bijections
$\ast,\ \omega: \mathfrak{U}\longleftrightarrow\mathfrak{V}$.
\item[(ii)]
There are bijections
$\ast,\ \omega:\ \mathfrak{U}_{[\geqslant\xi]}\longleftrightarrow \mathfrak{V}_{[\leqslant -\xi]}$ and
$\ast,\ \omega:\ \mathfrak{U}_{[>\xi]}\longleftrightarrow \mathfrak{V}_{[< -\xi]}$ for $\xi\in P^+$.
\end{itemize}
\end{lemma}
\begin{proof}
 To prove $(\romannumeral 1)$, it is sufficient to show that
$\ast(\widetilde{U}^{\prime})=\widetilde{U}^{\prime\prime}$
and $\omega(\widetilde{U}^{\prime})=\widetilde{U}^{\prime\prime}$.
For $b\in B^{\prime}$, $G(b)$ annihilates all $V(\lambda)$ for $\lambda\in P^+$. Then we have
$G(b)\circ {}^{\ast}V(-\lambda)^{\circ}=0$ for all $\lambda\in P^+$ which implies that
$G(b)^{\ast}=G(b^{\ast})$ annihilates all $V(-\lambda)$ for $\lambda\in P^+$. Hence
$\ast(\widetilde{U}^{\prime})\subseteq\widetilde{U}^{\prime\prime}$. Similarly we have
$\ast(\widetilde{U}^{\prime\prime})\subseteq\widetilde{U}^{\prime}$. It follows from
$\ast^2=id$ on $\widetilde{U}$ that $\ast(\widetilde{U}^{\prime})=\widetilde{U}^{\prime\prime}$. Given
$b\in B^{\prime}$, $G(b)$ annihilates all ${}^{\omega}V(-\lambda)$ for $\lambda\in P^+$. It implies that
$\omega(G(b))$ annihilates all $V(-\lambda)$ for $\lambda\in P^+$ and thus
$\ast(\widetilde{U}^{\prime})\subseteq\widetilde{U}^{\prime\prime}$. The proof of the equality is similar to that for $\ast$.
In order to prove $(\romannumeral 2)$, we only show $\ast(\mathfrak{U}_{[\geqslant\xi]})
\subseteq\mathfrak{V}_{[\leqslant -\xi]}$ and $\omega(\mathfrak{U}_{[\geqslant\xi]})
\subseteq\mathfrak{V}_{[\leqslant -\xi]}$. Given $x\in \mathfrak{U}_{[\geqslant\xi]}$, if $x^{\ast}$ acts on
$V(-\eta)$, $\eta\in P^+$, as a nonzero map, one can see that $x({}^{\ast}V(-\eta)^{\circ})=xV(\eta)\neq 0$ which implies
$\eta\geqslant\xi$. Hence $x^{\ast}\in \mathfrak{V}_{[\leqslant -\xi]}$. Similarly if $\omega(x)$ acts on $V(-\eta)$ for some $\eta\in P^+$,
as a nonzero map, then $\omega(x)V(-\eta)=x\circ{}^{\omega}V(-\eta)=xV(\eta)\neq 0$. Hence $\eta\geqslant \xi$ which implies
$\omega(x)\in \mathfrak{V}_{[\leqslant -\xi]}$.
\end{proof}

Note that $\mathfrak{U}=\bigoplus_{\xi\in P^+}\mathfrak{U}_{[\xi]}$ as a direct
sum of vector spaces. We have an isomorphism
\[\mathfrak{U}_{[\geqslant\xi]}/\mathfrak{U}_{[>\xi]}\cong
\mathfrak{U}_{[\xi]}\] as $k$-vector spaces. Furthermore,
$\mathfrak{U}_{[\geqslant\xi]}/\mathfrak{U}_{[>\xi]}$ is an algebra as well
as a $\mathfrak{U}$-bimodule which we call \emph{two-sided cell module} of
$\mathfrak{U}$ and denote also by $\mathfrak{U}(\xi)$ for simplicity.  This cell naturally inherits from $\mathfrak{U}$ a
canonical basis and its crystal basis is a family of copies of
$B(\xi)$. We have the following result similar to \cite{Lusztig:1993}.

\begin{proposition}
For $\xi\in P^+$,
\begin{itemize}
\item[(i)]
$\mathfrak{U}(\xi)$ decomposes into a direct sum of nice irreducible highest weight left $U$-submodules, each summand is isomorphic to
$V(\xi)$.
\item[(ii)]
$\mathfrak{U}(\xi)$ decomposes into a direct sum of nice irreducible highest weight right $U$-submodules, each summand is isomorphic to
$V(\xi)^{\circ}$.
\item[(iii)]
$\mathfrak{U}(\xi)\cong V(\xi)\otimes V(\xi)^{\circ}$ as $U$ or $\mathfrak{U}$-bimodules.
\end{itemize}
\end{proposition}

\begin{proof}
$(\romannumeral 1)$ is obvious. Since
we have bijections $\omega \circ \ast: \mathfrak{U}_{[\geqslant\xi]}\longleftrightarrow\mathfrak{U}_{[\geqslant\xi]}, \
 \mathfrak{U}_{[>\xi]}\longleftrightarrow\mathfrak{U}_{[>\xi]}$ by
Lemma~\ref{star}, $\omega\circ\ast$ induces an anti-automorphism of $\mathfrak{U}(\xi)$. Apply $\omega\circ\ast$ to
any summand $V$ in $(\romannumeral 1)$, we obtain a nice irreducible right $U$-module $\omega\circ\ast(V)$ by
Corollary~\ref{omega} and this proves $(\romannumeral 2)$.
Let $\phi$ be the restricting map on
$\mathfrak{U}_{[\geqslant\xi]}$ of the $\mathfrak{U}$-action on $V(\xi)$,
i.e. $\phi:\ \mathfrak{U}_{[\geqslant\xi]}\longrightarrow
End_{k}(V(\xi)).$ Then $\phi$ is a homomorphism of algebras without
1. It can be seen from Lemma~\ref{equvilant3} that the kernel of $\phi$ is
exactly $\mathfrak{U}_{[>\xi]}$ and thus $\bar{\phi}:\
\mathfrak{U}(\xi)\longrightarrow
End_{k}(V(\xi))$ is injective. We view $V(\xi)\otimes V(\xi)^{\circ}$ as a subset of $End_{k}(V(\xi))$, i.e.
\[(x\otimes f)(v)=f(v)x\ \textrm{for}\ x,\ v\in V(\xi),\ f\in V(\xi)^{\circ}. \]
It is easy to see that $V(\xi)\otimes V(\xi)^{\circ}$ is a $U$ or $\mathfrak{U}$-subbimodule as well as a subalgebra of
$End_{k}(V(\xi))$ where $U$ or $\mathfrak{U}$ acts on $V(\xi)\otimes V(\xi)^{\circ}$ as
\[(x(v\otimes f)y)(m)=f(ym)xv\ \textrm{for}\ v,m\in V(\xi),\ f\in V(\xi)^{\circ},\ x,y\in U\ \textrm{or}\ \mathfrak{U}.\]
 In fact the $\bar{\phi}$ defined above maps $\mathfrak{U}(\xi)$ injectively into $V(\xi)\otimes V(\xi)^{\circ}$, and moreover,
$\bar{\phi}:\ \mathfrak{U}(\xi)\longrightarrow V(\xi)\otimes V(\xi)^{\circ}$ is a homomorphism of $U$ or $\mathfrak{U}$-modules.
 Fixing a right weight $\eta\in P$,
$\mathfrak{U}(\xi)a_{\eta}$ is, by
Proposition~\ref{multiplicity},  a direct sum of $dimV(\xi)_{\eta}$ copies of $V(\xi)$
as a left $\mathfrak{U}$-module, where we denote the image of
$a_{\eta}$ in $\mathfrak{U}(\xi)$ by the same symbol. Hence $\bar{\phi}:\ \mathfrak{U}(\xi)\longrightarrow V(\xi)\otimes V(\xi)^{\circ}$ is surjective and
$\mathfrak{U}(\xi)\cong V(\xi)\otimes V(\xi)^{\circ}$.
\end{proof}

\section{Quantum coordinate ring}

\subsection{$\mathfrak{U}^{\circ}$}

For $U_1, U_2\in\{\widetilde{U},\ \widetilde{U}^{\prime},\ \mathfrak{U},\
U\textrm{-modules\ with\ canonical\ bases}\}$, let
$U_1\widehat{\otimes}U_1 $ be the set of all formal (possibly
infinite) linear combinations $\sum k_{b_1,b_2}G(b_1)\otimes
G(b_2).$ For $\lambda, \lambda_1, \lambda_2\in P$, the coproduct on
$U$ induces the map
$\Delta_{\lambda,\lambda_1,\lambda_2}: \
Ua_{\lambda}\longrightarrow Ua_{\lambda_1}\otimes Ua_{\lambda_2}$
where $\Delta_{\lambda,\lambda_1,\lambda_2}$ is nonzero only if
$\lambda=\lambda_1+\lambda_2$. Set
\[\Delta=\sum_{\lambda, \lambda_1, \lambda_2\in P}\Delta_{\lambda,\lambda_1,\lambda_2}:\ \widetilde{U}\longrightarrow
\widetilde{U}\widehat{\otimes}\widetilde{U}.\] For $a, b, c\in
\widetilde{B}$, we define $\hat{m}^{b,c}_a\in k$ to satisfy that
$\Delta(G(a))=\sum_{b,c}\hat{m}^{b,c}_aG(b)\otimes G(c).$
Note that $\hat{m}^{b,c}_a$ here is actually in $\Bbb Z[q,q^{-1}]$
\cite{Lusztig:1993}. For $\lambda$, $\lambda_1$, $\lambda_2$, $\mu$, $\mu_1$,
$\mu_2\in P^+$ with $\lambda=\lambda_1+\lambda_2$ and
$\mu=\mu_1+\mu_2$, let $\tau_1$, $\tau_2$ be the $U$-map
\[\tau_1: V(\lambda)\longrightarrow V(\lambda_1)\otimes V(\lambda_2),
 \ \ \ u_{\lambda}\longmapsto u_{\lambda_1}\otimes u_{\lambda_2},\ \ \ \
 \ \ \ \ \ \ \ \]
\[\tau_2: V(-\mu)\longrightarrow V(-\mu_1)\otimes V(-\mu_2),
 \ u_{-\mu}\longmapsto u_{-\mu_1}\otimes u_{-\mu_2}.\]
Set $R_{\lambda_2,-\mu_1}$ to be the unique isomorphism of
$U$-modules (R-matrix)
\[R_{\lambda_2,-\mu_1}:\ V(\lambda_2)\otimes V(-\mu_1)\longrightarrow\ V(-\mu_1)\otimes V(\lambda_2)\]
such that $R_{\lambda_2,-\mu_1}(u_{\lambda_2}\otimes
u_{-\mu_1})=u_{-\mu_1}\otimes u_{\lambda_2}$. Let $\tau$ be the
composition of $\tau_1\otimes\tau_2$ and $1\otimes
R_{\lambda_2,-\mu_1}\otimes 1$, i.e.
\[\xymatrix{
  V(\lambda, -\mu) \ar[dr]_{\tau} \ar[r]^-{\tau_1\otimes\tau_2}
                & \ \ V(\lambda_1)\otimes V(\lambda_2)\otimes V(-\mu_1)\otimes V(-\mu_2) \ar[d]^{1\otimes R_{\lambda_2,-\mu_1}\otimes 1}  \\
                & V(\lambda_1, -\mu_1)\otimes V(\lambda_2, -\mu_2)           }\]
Let $\rho$ be the map
$\rho: \ \widetilde{U}\widehat{\otimes}\widetilde{U}\longrightarrow V(\lambda_1,-\mu_1)
\widehat{\otimes} V(\lambda_2,-\mu_2)$ such that \[\rho(\sum
k_{a,b}G(a)\otimes G(b))=\sum
k_{a,b}(G(a)u_{\lambda_1,-\mu_1})\otimes
(G(b)u_{\lambda_2,-\mu_2}).\] One can see that $\widetilde{U}$ acts on
$u_{\lambda_1,-\mu_1}\otimes u_{\lambda_2,-\mu_2}$ as a map which
can be obtained through $\Delta$. More precisely, we have a
commutative diagram
\[\begin{CD}
\widetilde{U}
@>\Delta>> \widetilde{U}\widehat{\otimes}\widetilde{U} \\
@VV{\gamma}V @VV{\rho}V\\
V(\lambda_1,-\mu_1)\otimes V(\lambda_2,-\mu_2) @>i>>
V(\lambda_1,-\mu_1) \widehat{\otimes} V(\lambda_2,-\mu_2)
\end{CD}\]
where $\gamma(x)=x(u_{\lambda_1,-\mu_1}\otimes
u_{\lambda_2,-\mu_2})$ and $i$ is the canonical inclusion.

\begin{proposition}\label{commutative}
The following diagram is commutative
\[\begin{CD}
\widetilde{U}
@>\Delta>> Im(\Delta) \\
@VV{\alpha_{\lambda,-\mu}}V @VV{\rho|_{Im_{\Delta}}}V\\
V(\lambda,-\mu) @>\tau>> V(\lambda_1,-\mu_1) \otimes
V(\lambda_2,-\mu_2)
\end{CD} \ \ \ \ \ (5.1)\]
\end{proposition}

\begin{proof}
We regard $\widetilde{U}\widehat{\otimes}\widetilde{U}$ as a left
$\widetilde{U}$-module through $\Delta$. Hence all the maps are
homomorphisms of $\widetilde{U}$-modules. It is easy to check that
the two compositions in the diagram coincide when applied to
$a_{\xi}$ for any $\xi\in P$.
\end{proof}

\begin{lemma}\label{delta}
$\Delta(\widetilde{U}^{\prime})\subseteq
\widetilde{U}^{\prime}\widehat{\otimes}\widetilde{U}+\widetilde{U}\widehat{\otimes}\widetilde{U}^{\prime}.$
\end{lemma}

\begin{proof}
We assume that $\Delta(\widetilde{U}^{\prime})\nsubseteq
\widetilde{U}^{\prime}\widehat{\otimes}\widetilde{U}+\widetilde{U}\widehat{\otimes}\widetilde{U}^{\prime}.$
Then there exist $a\in B^{\prime}$ and $b, c\in
\widetilde{B}\setminus B^{\prime}$ such that $\hat{m}^{b,c}_a\neq
0$. We suppose that $b\in B(\lambda_1,-\mu_1)$ and $c\in
B(\lambda_2,-\mu_2)$ for some $\lambda_1$, $\lambda_2$, $\mu_1$,
$\mu_2\in P^+$. By Proposition~\ref{commutative}, we have
\[\rho\Delta(G(a))=\tau\alpha_{\lambda,-\mu}(G(a))\neq 0\] where
$\lambda=\lambda_1+\lambda_2$, $\mu=\mu_1+\mu_2$. Hence $a\in
B(\lambda,-\mu)$. It is known that there is a nice filtration of
$V(\lambda_i,-\mu_i)$
\[V(\lambda_i, -\mu_i)=F_{i,0}\supseteq F_{i,1}\supseteq\cdots\supseteq
F_{i,l}\supseteq\cdots\ \ \ \ \ \ \ \ \ \ \ \ \ (5.2(i))\] for $i=1,2$. Since $b,
c\notin B^{\prime}$, there exist $s$ and $t$ such that $G(b)\in
F_{1,s}$, $G(b)\notin F_{1,s+1}$ and $G(c)\in F_{2,t}$,
$G(c)\notin F_{2,t+1}.$ Let $\pi$ be the canonical map
\[\pi: \ V(\lambda_1,-\mu_1) \otimes
V(\lambda_2,-\mu_2)\longrightarrow (V(\lambda_1,-\mu_1)/F_{1,s+1})
\otimes (V(\lambda_2,-\mu_2)/F_{2,t+1})\] where
$(V(\lambda_1,-\mu_1)/F_{1,s+1}) \otimes
(V(\lambda_2,-\mu_2)/F_{2,t+1})$ is an object in $\mathcal O^+$,
denoted by $M$. Let $m$ be the image of $u_{\lambda_1,-\mu}\otimes
u_{\lambda_2,-\mu}$ in $M$. It follows from the assumption that
$\pi\tau\alpha_{\lambda,-\mu}(G(a))=G(a)m\neq 0.$ This contradicts
that $a\in B^{\prime}$ by Theorem~\ref{equvilant}.
\end{proof}

It follows from the lemma that $\Delta$ induces the map
$\bar{\Delta}: \ \mathfrak{U}\longrightarrow \mathfrak{U}\widehat{\otimes}\mathfrak{U}.$
For $a, b, c\in \mathfrak{B}$, similarly we define $\tilde{m}^{b,c}_a$ to
satisfy that
\[\bar{\Delta}(G(a))=\sum_{b,c}\tilde{m}^{b,c}_aG(b)\otimes G(c).\]
It is clear that $\hat{m}^{b,c}_a= \tilde{m}^{b,c}_a$ if $a, b, c\in
\mathfrak{B}=\widetilde{B}\setminus B^{\prime}$. The coassociation of the
coproduct on $U$ implies that
\[\sum_{c\in\mathfrak{B}}\tilde{m}^{a,b}_c\tilde{m}^{c,d}_e=\sum_{c\in\mathfrak{B}}\tilde{m}^{a,c}_e\tilde{m}^{b,d}_c
\ \ \ \ \ \ \ \ \ \ \ \ \ (5.3)\] for any $a, b, c, d, e\in\mathfrak{B}$.\\

One can see from the proof of Lemma~\ref{delta} that for $b\in
B(\lambda_1,-\mu_1)\bigcap\mathfrak{B}$ and $c\in
B(\lambda_2,-\mu_2)\bigcap\mathfrak{B}$ with $\tilde{m}^{b,c}_a\neq 0$,
we have $a\in B(\lambda_1+\lambda_2,-\mu_1-\mu_2)$. Hence when
$\mathfrak g$ is of finite type, the set $\{a\in\mathfrak{B}\ |\
\tilde{m}^{b,c}_a\neq 0\}$ is finite for fixed $b,c$ above since
$B(\lambda_1+\lambda_2,-\mu_1-\mu_2)$ is a finite set. We claim that it is
also true for $\mathfrak g$ of any type though
$B(\lambda_1+\lambda_2,-\mu_1-\mu_2)$ is not finite any more in
other cases.

\begin{theorem}\label{finiteness}
For $b, c\in\mathfrak{B}$, $\{a\in\mathfrak{B}\ |\ \tilde{m}^{b,c}_a\neq 0\}$
is a finite set.
\end{theorem}

\begin{proof} Assume that $b\in B(\lambda_1,-\mu_1)_{\xi_1}$ and $c\in
B(\lambda_2,-\mu_2)_{\xi_2}$ with $\tilde{m}^{b,c}_a\neq 0$, then
$a\in B(\lambda_1+\lambda_2,-\mu_1-\mu_2)_{\xi_1+\xi_2}$ from the
above statement. We suppose that $G(a)\in\mathfrak{U}_{[\xi]}$. As in the
proof of Lemma~\ref{delta}, there exist $s\geqslant 0$ and
$t\geqslant 0$ such that
$G(b)\in F_{1,s},\ G(b)\notin F_{1,s+1}$ and $G(c)\in
F_{2,t},\ G(c)\notin F_{2,t+1}$ where $F_{1,s}$, $F_{1,s+1}$ and
$F_{2,t}$, $F_{2,t+1}$ are in the filtration (5.2(1)) and (5.2(2))
respectively. Define $\pi$, $M$ and $m\in M$ as before. Hence we
obtain a homomorphism of $U$-modules
$\pi\tau:\ V(\lambda,-\mu)\longrightarrow M$ where
$\lambda=\lambda_1+\lambda_2$ and $\mu=\mu_1+\mu_2$. Since
$\tilde{m}^{b,c}_a\neq 0$, it implies that
\[\pi\tau(G(a)u_{\lambda,-\mu})=G(a)m\neq 0.\] Hence $\xi\leqslant\eta$
for some $\eta\in P^+$ such that $M_{\eta}\neq 0$ by
Lemma~\ref{equvilant2}. It follows that there exists an $l>0$ such
that
\[\{G(a)\ |\ \tilde{m}^{b,c}_a\neq 0\}\subseteq F_0\setminus F_{l}\]
where $F_0$, $F_l$ is in the filtration (3.2). Furthermore, there is
a bijection between $\{G(a)\ |\ \tilde{m}^{b,c}_a\neq 0\}$ and its
image in $F_0/F_l$ under the canonical map $\pi^{\prime}:
F_0\longrightarrow F_0/F_l$. Since $F_0/F_l\in ob(\mathcal O^+)$ and
\[\pi^{\prime}\{G(a)\ |\ \tilde{m}^{b,c}_a\neq 0\}\subseteq
(F_0/F_l)_{\xi_1+\xi_2},\] $\{G(a)\in \mathfrak{U}\ |\
\tilde{m}^{b,c}_a\neq 0\}$ is a finite set which proves the theorem.
\end{proof}

Let $\mathfrak{U}^{\ast}$ be the dual space of $\mathfrak{U}$, that is,
the set of all linear functions $\phi:\ \mathfrak{U}\longrightarrow k$.
For $b\in\mathfrak{U}$, set $b^{\ast}$ to be the linear function dual
to the canonical base element $G(b)$, i.e.
\[b^{\ast}(G(c))=\delta_{b,c}\ \textrm{for}\ b,c\in \mathfrak{B}.\] Let $\mathfrak{U}^{\circ}$ be the
subspace of $\mathfrak{U}^{\ast}$ spanned over $k$ by
$\{b^{\ast}\in\mathfrak{U}^{\ast}\ |\ b\in\mathfrak{B}\}$. We define an
algebra structure on $\mathfrak{U}^{\circ}$ by setting
\[b^{\ast}\cdot c^{\ast}=\sum_{a\in\mathfrak{B}}\tilde{m}^{b,c}_aa^{\ast}.\]
The sum is well-defined by Theorem~\ref{finiteness} and the
associativity of
this multiplication is implied by (5.3).\\

\subsection{Other versions of definition}

Let $U^{\ast}$ be the dual space of $U$, that is, the set of all
linear functions on $U$. The coproduct on $U$ provides a
multiplication on $U^{\ast}$, i.e. $(f_1\cdot f_2)(u)=\sum
f_1(u_{(1)})f_2(u_{(2)})$ where $\Delta(u)=\sum u_{(1)}\otimes
u_{(2)}$. For a $U$-module $M\in ob(\mathcal O^+)$, let $M^{\circ}$
be the graded dual of $M$ as before. For $m\in M$ and $f\in
M^{\circ}$, we define a \emph{coordinate function} $m\otimes f$ on
$U$ as
\[(m\otimes f)(u)=f(um)\ \textrm{for}\ u\in U.\]
\begin{definition}
The subalgebra of $U^{\ast}$ generated by all coordinate functions
$m\otimes f\in M\otimes M^{\circ}$ for all $M\in ob(\mathcal O^+)$
is called the quantum coordinate ring, denoted by $C_1$.
\end{definition}

For $m_i\otimes f_i\in M_i\otimes M_i^{\circ}$ with $M_i\in
ob(\mathcal O^+)$, $i=1, 2$, since \[(m_1\otimes
f_1)\cdot(m_2\otimes f_2)=(m_1\otimes m_2)\otimes(f_1\otimes f_2)\in
(M_1\otimes M_2)\otimes(M_1\otimes M_2)^{\circ}\] where $M_1\otimes
M_2\in ob(\mathcal O^+)$, $C_1$ is actually spanned by all
coordinate functions and has a structure of $U$-bimodule. The
complete irreducibility of category $\mathcal O^+$ implies the
following analogue of Peter-Weyl theorem.

\begin{proposition} (\cite{Kashiwara2:1993})
$C_1\cong \bigoplus_{\lambda\in P^+}V(\lambda)\otimes
V(\lambda)^{\circ}$ as $U$-bimodules and algebras.
\end{proposition}

Given $\lambda, \mu\in P^+$, we define a surjective map
$\tilde{\alpha}_{\lambda,-\mu}:\ U\longrightarrow M(\lambda,
-\mu)$ which takes $x$ to $xv_{\lambda, -\mu}$. Here
$M(\lambda,-\mu)$ inherits from $V(\lambda,-\mu)$ a canonical basis
which also has the stability property as in \cite{Lusztig:1992}, i.e. for
$\lambda, \mu, \theta\in P^+$, the $U$ (or $\mathfrak{U}$)-map
\[\tau_{\lambda+\theta,-\theta-\mu,\lambda,-\mu}:\
M(\lambda+\theta, -\theta-\mu)\longrightarrow M(\lambda, -\mu)\]
takes the canonical base elements to canonical base elements or
to zero. Note that $M(\lambda, -\mu)$ is usually not in $\mathcal O^+$
for $\mathfrak g$ of affine or indefinite type and more precisely,
its weight space might be infinite dimensional. We define
$M(\lambda, -\mu)^{\circ}$ to be the subspace of $M(\lambda,
-\mu)^{\ast}$ spanned by the dual basis associated to the canonical
basis of $M(\lambda, -\mu)$, i.e.
\[M(\lambda,-\mu)^{\circ}=\bigoplus_{b\in B(M(\lambda,-\mu))}kb^{\circ}\] where
$b^{\circ}(G(c))=\delta_{b,c}$ for $b,c\in B(M(\lambda,-\mu))$. One
can see that $\tilde{\alpha}_{\lambda,-\mu}$ induces the injective
map
\[\tilde{\alpha}_{\lambda,-\mu}^{\ast}:\ M(\lambda,-\mu)^{\circ}\longrightarrow
U^{\ast}.\] Denote by $U(\lambda,-\mu)^{\ast}$ the image of
$\tilde{\alpha}_{\lambda,-\mu}^{\ast}$. It follows from the
stability property of canonical bases that
\[U(\lambda,-\mu)^{\ast}\subseteq
U(\lambda+\theta,-\theta-\mu)^{\ast}.\] Indeed, for $b\in
B(M(\lambda,-\mu))\subseteq B(M(\lambda+\theta, -\theta-\mu))$, we
have
$\tilde{\alpha}_{\lambda,-\mu}^{\ast}(b^{\circ})=
\tilde{\alpha}_{\lambda+\theta,-\theta-\mu}^{\ast}(b^{\circ}),$
still denoted by $b^{\circ}$ in $U^{\ast}$ with no confusion.
\begin{definition}
$C_2$ is defined to be the subalgebra of $U^{\ast}$ generated by
$U(\lambda,-\mu)^{\ast}$ for all $\lambda,\mu\in P^+$.
\end{definition}

 Note that the above definition is a generalization of
Lusztig's definition of quantum coordinate rings for $\mathfrak g$
of finite type \cite{Lusztig:1992}.

\begin{lemma}
\begin{itemize}
\item[(i)]
$\tau : V(\lambda, -\mu)\longrightarrow V(\lambda_1, -\mu_1)\otimes
V(\lambda_2, -\mu_2)$ induces a homomorphism of $U$ or
$\mathfrak{U}$-modules
$\bar{\tau}:\ M(\lambda, -\mu)\longrightarrow M(\lambda_1,-\mu_1)\otimes M(\lambda_2,-\mu_2).$
\item[(ii)]
Let $\bar{\tau}^{\ast}$ be the map $\bar{\tau}^{\ast}:
(M(\lambda_1,-\mu_1)\otimes
M(\lambda_2,-\mu_2))^{\ast}\longrightarrow M(\lambda, -\mu)^{\ast}$
defined by $\bar{\tau}^{\ast}(f)(m)=f(\bar{\tau}(m))$ for $f\in
(M(\lambda_1,-\mu_1)\otimes M(\lambda_2,-\mu_2))^{\ast}$, $m\in
M(\lambda, -\mu)$. Then
$\bar{\tau}^{\ast}(M(\lambda_1,-\mu_1)^{\circ}\otimes M(\lambda_2,
-\mu_2)^{\circ})\subseteq M(\lambda, -\mu)^{\circ}$.
\end{itemize}
\end{lemma}

\begin{proof}
Let $\pi$ be the canonical $U$-map $\pi:\
V(\lambda_1,-\mu_1)\otimes V(\lambda_2,-\mu_2)\longrightarrow
M(\lambda_1,-\mu_1) \otimes M(\lambda_2,-\mu_2)$. Thus we have the composed
$U$-map $\pi\tau:\ V(\lambda,-\mu)\longrightarrow
M(\lambda_1,-\mu_1) \otimes M(\lambda_2,-\mu_2).$ Both
$M(\lambda_1,-\mu_1)$ and $M(\lambda_2,-\mu_2)$ have nice filtration
as in (4.1),
\[M(\lambda_i, -\mu_i)=M_{i,0}\supseteq M_{i,1}\supseteq\cdots\supseteq
M_{i,l}\supseteq\cdots\ \ \ \ \ \ \ \ \ \ \ \ \ \  (5.2(i))\] where $i=1,2$. For any $s, t\geqslant 0$, we define
$(\pi\tau)_{s,t}$ to be the composition of $\pi\tau$ with the canonical map
\[\eta_{s,t}:\ M(\lambda_1,-\mu_1)
\otimes M(\lambda_2,-\mu_2)\longrightarrow (M(\lambda_1,-\mu_1)
/M_{1,s})\otimes (M(\lambda_2,-\mu_2)/M_{2,t}),\] i.e.
$(\pi\tau)_{s,t}:\
V(\lambda,-\mu)\longrightarrow (M(\lambda_1,-\mu_1)/M_{1,s})
\otimes (M(\lambda_2,-\mu_2)/M_{2,t})$ where the right side is in $\mathcal O^+$.
Thus for any $b\in B^{\prime}(\lambda,-\mu)$, $(\pi\tau)_{s,t}(G(b))=0$.
It implies $(\pi\tau)(G(b))\in M_{1,s}\otimes M(\lambda_2,-\mu_2)+M(\lambda_1,-\mu_1)\otimes M_{2,t}$ for any $s, t\geqslant 0.$
Hence \[(\pi\tau)(G(b))\in\bigcap_{s,t\geqslant 0} (M_{1,s}\otimes M(\lambda_2,-\mu_2)+M(\lambda_1,-\mu_1)\otimes M_{2,t})=0.\]
It follows that $\pi\tau$ factors through $\bar{\tau}:\
M(\lambda, -\mu)\longrightarrow M(\lambda_1,-\mu_1)\otimes M(\lambda_2,-\mu_2)$ which proves $(\romannumeral 1)$.
To prove $(\romannumeral 2)$, we only show that
$\bar{\tau}^{\ast}(b_1^{\circ}\otimes b_2^{\circ})\subseteq M(\lambda,-\mu)^{\circ}$ for any
$b_i\in B(M(\lambda_i, -\mu_i))$, $i=1, 2$. Assume that
$\bar{\tau}^{\ast}(b_1^{\circ}\otimes b_2^{\circ})(G(b))\neq 0$ for some $b\in B(M(\lambda,-\mu))$,
that is, $(b_1^{\circ}\otimes b_2^{\circ})(\bar{\tau}G(b))\neq 0$. Then it
implies that $\tilde{m}_b^{b_1,b_2}\neq 0$. We know from
Theorem~\ref{finiteness} that there are only finitely many such
$G(b)$. Hence $\bar{\tau}^{\ast}(b_1^{\circ}\otimes
b_2^{\circ})\subseteq M(\lambda,-\mu)^{\circ}$.
\end{proof}

As a consequence, we have the following corollary.
\begin{corollary}\label{C_2}
\begin{itemize}
\item[(i)]
$U(\lambda_1,-\mu_1)^{\ast}\cdot U(\lambda_2,-\mu_2)^{\ast}\subseteq
U(\lambda_1+\lambda_2,-\mu_1-\mu_2)^{\ast}$.
\item[(ii)]
$C_2=\sum_{\lambda, -\mu\in P^+}U(\lambda,-\mu)^{\ast}$ and
$\{b^{\circ}\ |\ b\in \mathfrak{B} \}$ forms a basis of $C_2$.
\end{itemize}
\end{corollary}

\subsection{Equivalence of definitions}

For $x, y\in U$, $f\in \mathfrak{U}^{\circ}$ and $f\in \mathfrak{U}^{\circ}$,
we define $x\cdot f\cdot y\in \mathfrak{U}^{\ast}$ to satisfy that
\[(x\cdot f\cdot y)(u)=f(yux)\] for any $u\in \mathfrak{U}$. Suppose that
$x\in U_{\xi_1}$, $y\in U_{\xi_2}$ and $f=b^{\ast}$ such that
$G(b)\in \mathfrak{U}_{[\eta]}$ is in the image of $Ua_{\lambda}$ with
weight $\mu$. For $u=G(b^{\prime})\in\mathfrak{U}$, it can be seen from
the weight that $(x\cdot f\cdot y)(u)\neq 0$ implies $u$ is in the
image of $Ua_{\lambda+\xi_1}$ with weight $\mu-\xi_1-\xi_2$. Also we
have, by Corollary~\ref{ideal}, that
$u\in\oplus_{\theta\leqslant\eta}\mathfrak{U}_{[\theta]}$ if $(x\cdot
f\cdot y)(u)\neq 0$. It follows as in the proof of
Theorem~\ref{finiteness} that $x\cdot f\cdot y$ acts as zero for all
but finitely many $G(b^{\prime})\in \mathfrak{U}$. Hence $x\cdot f\cdot
y\in \mathfrak{U}^{\circ}$ and one can view $\mathfrak{U}^{\circ}$ as a
$U$-bimodule or similarly as a $\widetilde{U}$ or
$\mathfrak{U}$-bimodule. Fixing a left weight $\mu$, that is, taking
$x=a_{\mu}\in \widetilde{U}$, we can see from
Proposition~\ref{multiplicity} that the right
$\widetilde{U}$-submodule $a_{\mu}\mathfrak{U}^{\circ}$ of
$\mathfrak{U}^{\circ}$
 corresponds to a crystal which consists of $dimV(\lambda)_{\mu}$ copies of $B(\lambda)^{\circ}$  associated to irreducible integrable highest weight right $\widetilde{U}$-module $V(\lambda)^{\circ}$
for all $\lambda\in P^+$. Obviously the same happens to $C_1$ when applying $a_{\mu}$ to
the left side. \\

On the other hand, it is easy to see that the structure constant of multiplication
 in $C_2$ with respect to the basis $\{b^{\circ}\ |\ b\in
\mathfrak{B}\}$ is exactly the same as that in $\mathfrak{U}^{\circ}$ since
both multiplications are defined through coproduct on $U$. Thus
$\mathfrak{U}^{\circ}$ and $C_2$ are isomorphic as algebras. All the
statements above lead us to a belief that the three definitions of
quantum coordinate ring are equivalent.

\begin{theorem}\label{isomorphism}
$\mathfrak{U}^{\circ}\cong C_1\cong C_2$ as algebras.
\end{theorem}

\begin{proof}
We only show $C_1\cong C_2$. Given $f\in U(\lambda,
-\mu)^{\ast}\subseteq C_2$ with $\lambda, \mu\in P^+$, there is a
$g\in M(\lambda, -\mu)^{\circ}$ such that $f(x)=g(xv_{\lambda,
-\mu})$ for all $x\in U$. Since $g$ acts as zero for all but
finitely many canonical base elements of $M(\lambda, -\mu)$,
$g(M_s)=0$ for some $s\geqslant 0$ where $M_s$ is in the filtration
(4.1). Hence $g$ induces the linear map $\bar{g}:\
M(\lambda,-\mu)/M_s\longrightarrow k$ where $M(\lambda,-\mu)/M_s$,
denoted by $M$, is an object in $\mathcal O^+$. Clearly $\bar{g}\in
M^{\circ}$ and we have
\[f(x)=g(xv_{\lambda, -\mu})=\bar{g}(x(v_{\lambda,-\mu}+M_s))=((v_{\lambda,-\mu}+M_s)\otimes\bar{g})(x).\] We denote $v_{\lambda,-\mu}+M_s\in M$ by $m$. Thus $f=m\otimes \bar{g}\in M\otimes M^{\circ}$ which implies that $f\in C_1$. Conversely,
assume that $f=m\otimes g\in M\otimes M^{\circ}\subseteq C_1$ for some $M\in \mathcal O^+$ and $m\in M_{\xi}$.
There exists $\lambda, \mu\in P^+$ with $\lambda-\mu=\xi$ such that
$\phi:\ M(\lambda, -\mu)\longrightarrow M$ which takes $xv_{\lambda,-\mu}$ to $xm$ for $x\in U$ is a well-defined $U$-map.
Note that $\phi$ induces an injective map
\[\phi^{\ast}:\ (Im\phi)^{\ast}\longrightarrow M(\lambda,-\mu)^{\ast}.\]
As above we have $\phi(M_s)$ for some $s\geqslant 0$ where $M_s$ is in the filtration (4.1) and $\phi$ factors through a surjective
$U$-map $\bar{\phi}:\ M(\lambda,-\mu)/M_s\longrightarrow Im\phi$ which induces
\[\bar{\phi}^{\ast}:\ (Im\phi)^{\ast}\longrightarrow (M(\lambda,-\mu)/M_s)^{\ast}.\]
Since $Im\phi$, $M(\lambda,-\mu)/M_s\in \mathcal O^+$ and $g\in M^{\circ}$, then
$g|_{Im\phi}\in (Im\phi)^{\circ},\ \bar{\phi}^{\ast}(g|_{Im\phi})\in
(M(\lambda,-\mu)/M_s)^{\circ}.$ It implies that $\phi^{\ast}(g|_{Im\phi})\in
M(\lambda,-\mu)^{\circ}$. One can check easily that
\[f=m\otimes g=\tilde{\alpha}^{\ast}_{\lambda, -\mu}\phi^{\ast}
(g|_{Im\phi})\in C_2\] which proves the whole theorem.
\end{proof}
\begin{remark}
\begin{itemize}
\item[(i)]
 One may notice that the proof of Theorem~\ref{isomorphism} also implies the equality in Corollary~\ref{C_2} $(\romannumeral 2)$.
\item[(ii)]
Let $Z=\Bbb Z[q, q^{-1}]$. Similarly as that constructed by Lusztig
in \cite{Lusztig:2009} for $\mathfrak g$ of finite type, we can define a
$\Bbb Z$-form of the quantum coordinate ring by spanning a free
$Z$-module with the basis $\{b^{\ast}\ | b\in\mathfrak{B}\}$.
\end{itemize}
\end{remark}

\end{document}